\documentclass[11pt,draft]{article}
\usepackage{amsmath}
\usepackage{amssymb}

\newtheorem{thm}{Theorem}[section]
\newtheorem{defn}[thm]{Definition}
\newtheorem{lemma}[thm]{Lemma}
\newtheorem{prop}[thm]{Proposition}
\newtheorem{rem}[thm]{Remark}
\newtheorem{cor}[thm]{Corollary}
\begin{document}

\title{On the Laplace transform of absolutely monotonic functions}

\author{Stamatis Koumandos and Henrik L. Pedersen\footnote{Research
    supported by grant DFF–4181-00502
from The Danish Council for Independent Research $|$ Natural Sciences}}

\date{\today}
\maketitle

\begin{abstract}
We obtain necessary and sufficient conditions on a function in order that it be the Laplace transform of an absolutely monotonic function. Several closely related results are also given.
\end{abstract}


\section{Introduction and results}
Before stating our main results let us recall the notions of absolute and of complete monotonicity. A function $\varphi:[0,\infty)\to \mathbb R$ is called absolutely monotonic if it is infinitely often differentiable on $[0,\infty)$ and $\varphi^{(k)}(x)\geq 0$ for all $k\geq 0$ and all $x\geq 0$. It is well-known that an absolutely monotonic function $\varphi$ on $[0,\infty)$ has an extension to an entire function with the power series expansion 
$
\varphi(z)=\sum_{n=0}^{\infty}a_nz^n,
$
where $a_n\geq 0$ for all $n\geq 0$. The Laplace transform of $\varphi$ is defined exactly when $\varphi$ extends to an entire function of at most exponential type zero, meaning that $\varphi$ has the following property. For any given $\epsilon>0$ there exists a positive constant $C_{\epsilon}$ such that 
$
|\varphi(z)|\leq C_{\epsilon}e^{\epsilon |z|}
$
for all $z$. We remark that the Laplace transform of an entire function of exponential type is defined on a half-line and is called the Borel transform of the entire function. It is related to the so-called indicator diagram of a function of exponential type, see \cite{boas}.

A function $f:(0,\infty)\to \mathbb R$ is called completely monotonic if $f$ is infinitely often differentiable on $(0,\infty)$ and $(-1)^nf^{(n)}(x)\geq 0$ for all $n\geq 0$ and all $x>0$. Bernstein's theorem states that $f$ is completely monotonic if and only if there exists a positive measure $\mu$ on $[0,\infty)$ such that $t\mapsto e^{-xt}$ is integrable w.r.t.\ $\mu$ for all $x>0$ and
$$
f(x)=\int_0^{\infty}e^{-xt}\, d\mu(t).
$$
The class of Stieltjes functions is an important subclass of the completely monotonic functions: A function $f:(0,\infty)\to \mathbb R$ is a Stieltjes function if 
$$
f(x)=\int_0^{\infty}\frac{d\mu(t)}{x+t}+c,
$$
where $\mu$ is a positive measure on $[0,\infty)$ making the integral converge for $x>0$ and $c\geq 0$. It is well known that the class of Stieltjes functions can be described as the Laplace transforms of completely monotonic functions.
Widder characterized this class as follows: $f$ is a Stieltjes function if and only if the function $(x^kf(x))^{(k)}$ is completely monotonic for all $k\geq 0$. (See \cite{Wid1}, and also \cite{Sok}.)

In order to motivate our results let us consider the following example. 
It is easy to show that 
$
H(x)=x^{-1}e^{1/x}
$
satisfies
$$
H_k(x)\equiv (-1)^k\big(x^k\,H(x)\big)^{(k)}=x^{-(k+1)}\,e^{1/x}, \;\;x>0.
$$
Hence $H_k$ is completely monotonic for all $k\geq 0$, being a product of completely monotonic functions. We also have
$$
H(x)=\int_0^{\infty}e^{-xt}h(t)\, dt,
$$
where $h(t)=\sum_{n=0}^{\infty}(n!)^{-2}t^n$ is absolutely monotonic.
This  example indicates that there must be an analogue of Widder's characterization mentioned above for Laplace transforms of absolutely monotonic functions and we present it in Theorem \ref{thm:main1} below. In our analogue it will become clear that the situation is different than in Widders characterization.
\begin{thm}
\label{thm:main1}
The following properties of a function $f:(0,\infty)\to \mathbb R$ are equivalent.
\begin{enumerate}
\item[(i)] There is an absolutely monotonic function $\varphi:[0,\infty)\to \mathbb R$ such that
$$
f(x)=\mathcal L(\varphi)(x)=\int_0^{\infty}e^{-xt}\varphi(t)\, dt, \quad x>0.
$$
\item[(ii)] There is a sequence $\{a_n\}$, with $a_n\geq 0$, such that we have for all $n\geq 0$ 
$$
f(x)=\sum_{k=1}^{n}\frac{a_k}{x^k}+R_n(x),\quad x>0
$$
where $R_n$ is a completely monotonic function of order $n$.
\item[(iii)] The function $(-1)^k(x^k\,f(x))^{(k)}$ is completely monotonic for all $k\geq 0$. 
\item[(iv)] The function $(-1)^k(x^k\,f(x))^{(k)}$ is non-negative for all $k\geq 0$.
\item[(v)] We have $f(x)\geq 0$ and $(x^k\,f(x))^{(2k-1)}\leq 0$ for all $k\geq 1$. 
\end{enumerate}
\end{thm}
If the Laplace transform $f$ of an absolutely monotonic function is defined on a half line then $f$ has the properties (ii), (iii), (iv) and (v) above. Moreover, (iii), (iv) and (v) are equivalent for a function $f$ defined on a half line $(a,\infty)$.
 
A function $f:(0,\infty)\to \mathbb R$ is called completely monotonic of order $\alpha>0$ if $x^{\alpha}f(x)$ is completely monotonic. These functions have been studied in different contexts and in particular as remainders in asymptotic formulae, see \cite{kp1} and \cite{kp2}.

To show that (iii) implies (i) in Theorem \ref{thm:main1} we characterize those functions $f$ that satisfy (iii) up to some given positive integer, in terms of the properties of the representing measure of $f$ itself. More specifically, we shall prove theorems \ref{thm:main2} and \ref{thm:main3} below. Before stating these results we recall some classes introduced in \cite{kp2}.
\begin{defn}
Let $A_0$ denote the set of positive Borel measures $\sigma$ on $[0,\infty)$ such that $\int_0^{\infty}e^{-xs}\,d\sigma(s)<\infty$ for all $x>0$, let $A_1(\sigma)$ denote the set of functions $t\mapsto \sigma([0,t])$, where $\sigma\in A_0$, and for $n\geq 2$, let $A_n(\sigma)$ denote the set of $n-2$ times differentiable functions $\xi:[0,\infty)\to \mathbb R$ satisfying $\xi^{(j)}(0)=0$ for $j\leq n-2$ and $\xi^{(n-2)}(t)= \int_0^t\sigma([0,s])\, ds$ for some $\sigma\in A_0$. 
\end{defn}
These classes are simply the fractional integrals of positive integer order of a Borel measure. 
\begin{thm}
\label{thm:main2}
Let $f:(0,\infty)\to \mathbb R$ and $k\geq 1$. If 
$(-1)^j(x^jf(x))^{(j)}$ is completely monotonic for $j=0,\ldots,k$ then 
$$
f(x)=\int_0^{\infty}e^{-xt}p(t)\, dt,
$$
where for some measures $\sigma,\sigma_0,\ldots,\sigma_{k-1}\in \mathcal A_0$ we have 
$t^{k-1}p(t)\in \mathcal A_k(\sigma)$, 
$t^{j}p^{(j)}(t)\in \mathcal A_1(\sigma_j)$ for $0\leq j\leq k-1$ and where 
the measure 
$$
d\sigma_k(t)=td\sigma_{k-1}(t)-(k-1)\sigma_{k-1}([0,t])dt
$$
is positive. Furthermore, 
\begin{align*}
(-1)^j(x^jf(x))^{(j)}&=\int_0^{\infty}e^{-xt}t^{j}p^{(j)}(t)\, dt,\quad \text{for}\ j\leq k-1,\ \text{and}\\
(-1)^k(x^kf(x))^{(k)}&=\int_0^{\infty}e^{-xt}\, d\sigma_k(t).
\end{align*}
\end{thm}
This theorem has a converse.
\begin{thm}
\label{thm:main3}
Let $k\geq 1$ be given. Suppose that 
$$
f(x)=\int_0^{\infty}e^{-xt}p(t)\, dt,
$$
where $t^jp^{(j)}(t)\in \mathcal A_1(\mu_j)$ for all $j\leq k-1$ and that 
$$
d\mu_k(t)=t\, d\mu_{k-1}(t)-(k-1)\mu_{k-1}([0,t])dt
$$
is a positive measure. Then  $(-1)^j(x^jf(x))^{(j)}$ is completely monotonic for $j\leq k$, and 
\begin{align*}
(-1)^j(x^jf(x))^{(j)}&=\int_0^{\infty}e^{-xt}t^jp^{(j)}(t)\, dt,\quad j\leq k-1,\\
(-1)^k(x^kf(x))^{(k)}&=\int_0^{\infty}e^{-xt}\, d\mu_k(t).
\end{align*}
\end{thm}
Sokal \cite{Sok} introduced for $\lambda >0$ the operators
$$
T_{n,k}^{\lambda}(f)(x)\equiv (-1)^nx^{-(n+\lambda-1)}\left(x^{k+n+\lambda-1}f^{(n)}(x)\right)^{(k)}, \quad n,k\geq 0
$$
and showed that $f$ is of the form 
$$
f(x)=\int_0^{\infty}e^{-xt}t^{\lambda-1}\varphi(t)\, dt, \quad x>0
$$
for some completely monotonic function $\varphi$ if and only if $T_{n,k}^{\lambda}(f)(x)\geq 0$ for all $x>0$, and  $n,k\geq 0$. It turns out that the corresponding result when $\varphi$ is absolutely monotonic is relatively simple. We give it in Theorem \ref{thm:main4}.
\begin{thm}
\label{thm:main4}
Let $\lambda >0$ be given. The following properties of a function $f:(0,\infty)\to \mathbb R$ are equivalent.
\begin{enumerate}
\item[(i)] There is an absolutely monotonic function $\varphi:[0,\infty)\to \mathbb R$ such that
$$
f(x)=\int_0^{\infty}e^{-xt}t^{\lambda-1}\varphi(t)\, dt, \quad x>0.
$$
\item[(ii)] The function $(-1)^k(x^{k+\lambda-1}\,f(x))^{(k)}$ is completely monotonic for all $k\geq 0$. 
\item[(iii)] The function $(-1)^k(x^{k+\lambda-1}\,f(x))^{(k)}$ is non-negative for all $k\geq 0$.
\end{enumerate}
\end{thm}
The proof of this theorem follows from Theorem \ref{thm:main1} by noticing
$$
f(x)=\int_0^{\infty}e^{-xt}t^{\lambda -1}\varphi(t)\, dt
$$
for some absolutely monotonic function $\varphi$ of exponential type zero if and only if
$$
x^{\lambda -1}f(x)=\int_0^{\infty}e^{-xt}\psi(t)\, dt
$$
for some absolutely monotonic function $\psi$ of exponential type zero. Indeed,  the relationship between the functions $\varphi$ and $\psi$ is:
$$
\varphi(t)=\sum_{n=0}^{\infty}a_nt^n\Leftrightarrow \psi(t)=\sum_{n=0}^{\infty}\frac{a_n\Gamma(n+\lambda)}{n!}t^n.
$$

\section{Proof of Theorem \ref{thm:main1}}

{\it Proof that (i) implies (ii):} We have $\varphi(t)=\sum_{k=0}^{\infty}a_kt^k$, with $a_k\geq 0$ and get by interchanging integration and summation 
$$
f(x)=\int_{0}^{\infty}e^{-xt}\varphi(t)\,dt=\sum_{k=0}^{\infty}\frac{a_kk!}{x^{k+1}}=\sum_{k=0}^{n-1}\frac{a_kk!}{x^{k+1}}+\frac{1}{x^n}\sum_{k=n}^{\infty}\frac{a_kk!}{x^{k-n}},
$$
where the right-most sum is completely monotonic of order $n$.\hfill $\square$


{\it Proof that (ii) implies (iii):} This follows from the relation 
$$
(-1)^k(x^kf(x))^{(k)}=(-1)^k(x^kR_k(x))^{(k)},
$$ 
since $x^kR_k(x)$, and hence $(-1)^k(x^kR_k(x))^{(k)}$ is completely monotonic.\hfill $\square$


{\it Proof that (iii) implies (i):}  Let $k\geq 2$  be given. By Theorem \ref{thm:main2} we have 
$$
f(x)=\int_0^{\infty}e^{-xt}q_k(t)\, dt,
$$
where $q_k$ is $k-2$ times differentiable on $(0,\infty)$ with $q_k^{(j)}(t)\geq 0$ for all $j\leq k-2$ and all $t>0$.  Furthermore, $q_k$ is continuous on $[0,\infty)$, and it does not depend on $k$. \hfill $\square$

{\it Proof that (iii) implies (v):} The identity 
$$
-(x^kf(x))^{(2k-1)}=(-1)^{k-1}((-1)^k(x^kf(x))^{(k)})^{(k-1)}
$$
shows that $-(x^kf(x))^{(2k-1)}$ is completely monotonic, hence non-negative.\hfill $\square$

Let us introduce some notation and give two simple lemmas. For a sufficiently smooth function $f:(0,\infty)\to \mathbb R$ we define
\begin{equation}
\label{eq:fj}
f_j(x)=(-1)^j(x^jf(x))^{(j)}.
\end{equation}
\begin{lemma}
\label{lemma:hlp1}
We have $f_j(x)=-xf_{j-1}'(x)-jf_{j-1}(x)$, for $j\geq 1$.
\end{lemma}
\begin{lemma}
\label{lemma:hlp2}
We have $f_j^{(k-1)}(x)=-xf_{j-1}^{(k)}(x)-(j+k-1)f_{j-1}^{(k-1)}(x)$ for $j, k\geq 1$.
\end{lemma}
(Lemma \ref{lemma:hlp1} is easily proved using  $f_j(x)=-((x(-1)^{j-1}x^{j-1}f(x))^{(j-1)})'$ and Lemma \ref{lemma:hlp2} follows directly from Lemma \ref{lemma:hlp1} by differentiation.)

{\it Proof that (iv) implies (iii):} From Lemma \ref{lemma:hlp1} we get 
$$
-xf_k'(x)=f_{k+1}(x)+(k+1)f_k(x)\geq 0
$$ 
so $-f_k'(x)\geq 0$. If $(-1)^nf_k^{(n)}(x)\geq 0$ for some $n\geq 0$ and all $k\geq 0$ then
by Lemma \ref{lemma:hlp2} 
$$
x(-1)^{n+1}f_{k}^{(n+1)}(x)=(-1)^n\left( f_{k+1}^{(n)}(x)+(k+n+1)f_k^{(n)}(x)\right)\geq 0,
$$
and in this way it follows that $f_k$ is completely monotonic.\hfill $\square$

Before completing the proof of Theorem \ref{thm:main1} we need another lemma.
\begin{lemma}
\label{lemma:1}
Suppose that $f\,:\,(0,\,\infty) \rightarrow \mathbb{R}$  satisfies
$f(x)\geq 0$ and  
$$
(x^kf(x))^{(2k-1)}\leq 0 \quad \text{for all}\quad  k\geq 1.
$$
Then 
\begin{equation}\label{eq:6}
\lim_{x\to \infty}x^{k}f^{(k)}(x)=0,\quad \text{for all}\quad k\geq 0,
\end{equation}
and 
\begin{equation}
\label{eq:8}
\lim_{x\to \infty}\big(x^{k}\,f(x)\big)^{(\nu)}=0\quad \text{for all}\quad  \nu\geq k\geq 0.
\end{equation}
\end{lemma}
\emph{Proof:}
Since $f(x)\geq 0$ and $x\,f(x)$ is decreasing on $(0,\,\infty)$ we have that
$$
\lim_{x\to +\infty}x\,f(x)=B\geq 0.
$$
The following identity is given in \cite[Lemma 3.11]{Wid1}:
\begin{equation}\label{eq:7}
x^{k-1}\,\big(x^k\,f(x)\big)^{(2k-1)}=\big(x^{2k-1}\,f^{(k-1)}(x)\big)^{(k)}.
\end{equation}
Hence by assumption $(x^{2k-1}\,f^{(k-1)}(x))^{(k)}\leq 0$,
and therefore there exista a constant $c_{k}$ such that
\begin{equation*}
\big(x^{2k-1}\,f^{(k-1)}(x)\big)^{(k-1)}\leq c_{k}\;\;\text{for all}\;\; x\in[1,\,\infty)
\end{equation*}
Integrating this on $[1,\,x]$ we have that
\begin{equation*}
\big(x^{2k-1}\,f^{(k-1)}(x)\big)^{(k-2)}\leq O(x),\;\;\; x\to\infty
\end{equation*}
and repeating the same process we end up with
\begin{equation*}
x^{2k-1}\,f^{(k-1)}(x)\leq O(x^{k-1}),\;\;\; x\to\infty,
\end{equation*}
that is
\begin{equation*}
f^{(k-1)}(x)\leq O(x^{-k}),\;\;\; x\to\infty.
\end{equation*}
We next show by a Tauberian argument that these one-sided estimates yield \eqref{eq:6}. Put 
$g(x)=f(x)-Bx^{-1}$. Then $g(x)=o(x^{-1})$ and 
$$
g''(x)=f''(x)-2Bx^{-3}\leq O(x^{-3}),\quad \text{as}\ x\to \infty.
$$
By \cite[Theorem 4.4, page 193]{Wid2} we conclude 
$$
f'(x)+Bx^{-2}=g'(x)=o(x^{-2}),
$$
and in particular $xf'(x)\to 0$ as $x\to \infty$. This process is continued: $g'''(x)\leq O(x^{-4})$ and $g'(x)=o(x^{-2})$ gives $f''(x)-2Bx^{-3}=g''(x)=o(x^{-3})$, so that $x^2f''(x)\to 0$ for $x\to \infty$. In this way \eqref{eq:6} follows.
The second assertion follows from the first and Leibniz' rule.\hfill $\square$

{\it Proof that (v) implies (iv):} We show that $(-1)^k(x^kf(x))^{(k)}\geq 0$ for all $k\geq 0$. Let $k\geq 1$ be given. By (v), the function $(x^kf(x))^{(2k-2)}$ is decreasing on $(0,\infty)$ and by Lemma \ref{lemma:1} we see that $\lim_{x\to \infty} (x^kf(x))^{(2k-2)}=0$. Therefore, the function $(-1)^2(x^kf(x))^{(2k-2)}=(x^kf(x))^{(2k-2)}$ must be non-negative. This means in turn that the function $(x^kf(x))^{(2k-3)}$ is increasing on $(0,\infty)$ and by Lemma \ref{lemma:1} we see that $\lim_{x\to \infty} (x^kf(x))^{(2k-3)}=0$. Therefore, the function $(-1)^3(x^kf(x))^{(2k-3)}=-(x^kf(x))^{(2k-3)}$ must be non-negative. This argument is continued until we conclude that the function $(-1)^k(x^kf(x))^{(2k-k)}=(-1)^k(x^kf(x))^{(k)}$ must be non-negative.\hfill $\square$

\section{Proof of Theorem \ref{thm:main2} and Theorem \ref{thm:main3}}
Let us start by noticing some consequences of Lemma \ref{lemma:hlp1} and Lemma \ref{lemma:hlp2}, for the functions in \eqref{eq:fj}.
\begin{cor}
\label{cor:hlp1}
We have $(-xf_{j-1}(x))'=f_{j}(x)+(j-1)f_{j-1}(x)$ for $j\geq 1$.
\end{cor}
This follows immediately from Lemma \ref{lemma:hlp1}.
\begin{cor}
\label{cor:hlp2}
Let $k\geq 1$. For any $j\in \{0,\ldots,k-1\}$ we have 
$$
(-1)^jx^jf^{(k-1)}(x)=\sum_{l=0}^ja_{j,l}f_l^{(k-j-1)}(x),
$$
where $a_{j,l}\geq 0$.
\end{cor}
    {\it Proof:} For $j=0$ the assertion clearly holds, with $a_{0,0}=1$. Assume now that the assertion holds for $j$ and $j\leq k-2$. Then 
$$
(-1)^{j+1}x^{j+1}f^{(k-1)}(x)=-x\sum_{l=0}^ja_{j,l}f_l^{(k-j-1)}(x).
$$
Now, $-xf_l^{(k-1-j)}(x)=f_{l+1}^{(k-j-2)}(x)+(k-1-j+l)f_{l}^{(k-j-2)}(x)$ by Lemma \ref{lemma:hlp2} so
$$
(-1)^{j+1}x^{j+1}f^{(k-1)}(x)=\sum_{l=0}^ja_{j,l}(f_{l+1}^{(k-j-2)}(x)+(k-1-j+l)f_{l}^{(k-j-2)}(x))
$$
and the proof is complete.\hfill $\square$

\begin{lemma}
\label{lemma:hlp3}
Assume that $f_j$ is completely monotonic for all $j\in \{0,\ldots,k\}$. Then 
\begin{enumerate}
\item[(a)] $xf_j(x)$ is completely monotonic for $0\leq j\leq k-1$;
\item[(b)] $(-1)^{j-1}x^jf^{(j-1)}(x)$ is completely monotonic for $1\leq j\leq k$.
\end{enumerate}
\end{lemma}
{\it Proof:} To prove (a) we notice that the function $xf_j(x)$ is non-negative since $f_j$ is completely monotonic. By Corollary \ref{cor:hlp1},
$$
-(xf_j(x))'=f_{j+1}(x)+jf_j(x),
$$
so that $-(xf_j(x))'$ is completely monotonic.

To prove (b) we use Corollary \ref{cor:hlp2}:
$$
(-1)^{j-1}x^{j-1}f^{(j-1)}(x)=\sum_{l=0}^{j-1}a_{j-1,l}f_l(x),
$$
from which it follows that $(-1)^{j-1}x^{j-1}f^{(j-1)}(x)$ is completely monotonic. Furthermore,
$$
(-1)^{j-1}x^{j}f^{(j-1)}(x)=x\sum_{l=0}^{j-1}a_{j-1,l}f_l(x)
=\sum_{l=0}^{j-1}a_{j-1,l}xf_l(x),
$$
and for $l\leq j-1$, the function $xf_l(x)$ is completely monotonic by (a).\hfill $\square$

{\it Proof of Theorem \ref{thm:main2}:} By (b) of Lemma \ref{lemma:hlp3}, $(-1)^{k-1}f^{(k-1)}(x)$ is completely monotonic of order $k$, so 
$$
(-1)^{k-1}f^{(k-1)}(x)=\int_0^{\infty}e^{-xt}p_k(t)\, dt,
$$
where $p_k\in \mathcal A_k(\sigma)$. By the same lemma,
$xf(x)$ is also completely monotonic, so that $f$ is of the form
$$
f(x)=\int_0^{\infty}e^{-xt}p(t)\, dt
$$
for some $p\in \mathcal A_1(\sigma_0)$. Hence $(-1)^{k-1}f^{(k-1)}=\mathcal L(t^{k-1}p(t))$ so that $t^{k-1}p(t)=p_k(t)\in \mathcal A_k(\sigma)$.

Integration by parts yields
$$
f(x)=\int_0^{\infty}e^{-xt}p(t)\, dt=\frac{\sigma_0(\{0\})}{x}+\frac{1}{x}\int_0^{\infty}e^{-xt}\, d\sigma_0(t),
$$
which implies
$$
f_1(x)=(-xf(x))'=\int_0^{\infty}e^{-xt}t\, d\sigma_0(t).
$$
By Lemma \ref{lemma:hlp3} $f_1$ is completely monotonic of order 1, so that $t\, d\sigma_0(t)=\sigma_1([0,t])\, dt$. This means $tp'(t)\in \mathcal A_1(\sigma_1)$ and $f_1(x)=\int_0^{\infty}e^{-xt}tp'(t)\, dt$.

Assume now that we have obtained 
$$
t^jp^{(j)}(t)=\sigma_j([0,t])
$$
for some $j\leq k-2$, and $f_j(x)=\int_0^{\infty}e^{-xt}t^jp^{(j)}(t)\, dt$.
Integration by parts gives us
\begin{align*}
f_j(x)&=\frac{\sigma_j(\{0\})}{x}+\frac{1}{x}\int_0^{\infty}e^{-xt}\, d\sigma_j(t)\\
&=\frac{\sigma_j(\{0\})}{x}+\frac{1}{x}\int_0^{\infty}e^{-xt}\left(jt^{j-1}p^{(j)}(t)+t^jp^{(j+1)}(t)\right)\, dt,
\end{align*}
so that  
$$
(-xf_j(x))'= \int_0^{\infty}e^{-xt}t\,\left(jt^{j-1}p^{(j)}(t)+t^jp^{(j+1)}(t)\right)\, dt.
$$
But also $f_{j+1}$ is completely monotonic of order 1 so $f_{j+1}=\mathcal L(q)$, where $q\in \mathcal A_1(\sigma_{j+1})$ for some $\sigma_{j+1}$ and furthermore 
$$
(-xf_j(x))'= f_{j+1}(x)+jf_j(x).
$$
Comparing these two relations we find $t^{j+1}p^{(j+1)}(t)=q(t)\in \mathcal A_1(\sigma_{j+1})$. Finally,
\begin{align*}
f_k(x)&=(-xf_{k-1}(f))'-(k-1)f_{k-1}(x)\\
&=-\left(x\int_0^{\infty}e^{-xt}\sigma_{k-1}([0,t])\, dt\right)'-(k-1)\int_0^{\infty}e^{-xt}\sigma_{k-1}([0,t])\, dt\\
&=-\left(x\left(\frac{\sigma_{k-1}(\{0\})}{x}+\frac{1}{x}\int_0^{\infty}e^{-xt}\, d\sigma_{k-1}(t)\right)\right)'\\
&\phantom{=} \ -(k-1)\int_0^{\infty}e^{-xt}\sigma_{k-1}([0,t])\, dt\\
&=\int_0^{\infty}e^{-xt}\left\{td\sigma_{k-1}(t)-(k-1)\sigma_{k-1}([0,t])dt\right\}
\end{align*}
and the assumed complete monontonicity forces the representing measure to be positive.\hfill $\square$

{\it Proof of Theorem \ref{thm:main3}:}
The function $f_0=f$ is completely monotonic. Assume next that $f_j$ is completely monotonic with 
$$
f_j(x)=\int_0^{\infty}e^{-xt}t^jp^{(j)}(t)\, dt,
$$
for some $j\leq k-2$. Since 
\begin{align*}
f_{j+1}(x)&=-xf_{j}'(x)-(j+1)f_{j}(x)\\
&=x\int_0^{\infty}e^{-xt}t^{j+1}p^{(j)}(t)\, dt-(j+1)\int_0^{\infty}e^{-xt}t^{j}p^{(j)}(t)\, dt,
\end{align*}
and
\begin{align*}
\int_0^{\infty}e^{-xt}t^{j+1}p^{(j)}(t)\, dt&=\frac{1}{x}\int_0^{\infty}e^{-xt}t^{j+1}p^{(j+1)}(t)\, dt\\
&\phantom{=}\ +(j+1)\frac{1}{x}\int_0^{\infty}e^{-xt}t^{j}p^{(j)}(t)\, dt
\end{align*}
we obtain 
$$
f_{j+1}(x)=\int_0^{\infty}e^{-xt}t^{j+1}p^{(j+1)}(t)\, dt.
$$
Finally,
\begin{align*}
\int_0^{\infty}e^{-xt}t^{k}p^{(k-1)}(t)\, dt&=\int_0^{\infty}e^{-xt}t\mu_{k-1}([0,t])\, dt\\
&=\frac{1}{x}\int_0^{\infty}e^{-xt}\mu_{k-1}([0,t])\, dt+ \frac{1}{x}\int_0^{\infty}e^{-xt}t\,d\mu_{k-1}(t)
\end{align*}
and if we combine this with the relation $f_{k}(x)=-xf_{k-1}'(x)-kf_{k-1}(x)$ we get 
$$
f_k(x)=\int_0^{\infty}e^{-xt}\left\{td\mu_{k-1}(t)-(k-1)\mu_{k-1}([0,t])dt\right\},
$$
and hence also $f_k$ is completely monotonic.\hfill $\square$

\section{Concluding remarks}
\begin{rem}
Suppose that the function $f$ satisfies condition (iii) of Theorem \ref{thm:main1} and that there exists a positive integer $r\geq 2$ such that $\left(x^rf(x)\right)'\leq 0$ for all $x>0$. Then $f$ is completely monotonic of order $r$, and hence of the form
$$
f(x)=\int_0^{\infty}e^{-xt}\varphi(t)\, dt,
$$
where $\varphi$ is absolutely monotonic and satifies $\varphi^{(k)}(0)=0$ for all $k\leq r-2$.


{\it Proof:} We shall use the fact that a function $g$ is completely monotonic if $g(x)\geq 0$, $g'(x)\leq 0$, and $(-1)^mg^{(m)}(x)\geq 0$ for infinitely many $m$. See \cite[Corollary 1.14, p.~12]{SSV}. It then suffices to consider $g(x)=x^rf(x)$. The integral representation follows from (v) of Theorem \ref{thm:main1}. \hfill $\square$


\end{rem}
\begin{rem}
Suppose that $f$ is the Laplace transform of a polynomial $\varphi$ of degree $r-1$ with non-negative coefficients. Then clearly 
$$
f(x)=\sum_{k=0}^{r-1}\varphi^{(k)}(0)/x^{k+1}
$$ 
whence
$$
\lim_{x\to 0^+}x^{r}f(x)=\varphi^{(r-1)}(0)<\infty.
$$
Conversely, if the function $f$ satisfies any condition of Theorem \ref{thm:main1} and there exists a positive integer $r$ and a number $a_r\neq 0$ such that 
\begin{equation}
\label{eq:ar}
\lim_{x\to 0^+}x^{r}f(x)=a_r
\end{equation}
then $f$ is the Laplace transform of a polynomial of degree $r-1$ with non-negative coefficients. Indeed, we have 
$$
x^{r}f(x)=\int_0^{\infty}e^{-xt}\varphi^{(r)}(t)\, dt+\sum_{j=0}^{r-1}\varphi^{(r-j-1)}(0)x^{j},
$$
where $\varphi$ is an absolutely monotonic function. Suppose for a contradiction that $\varphi^{(r)}$ is not identically zero. Then, since $\varphi^{(r)}$ is increasing and hence not integrable,
$$
\lim_{x\to 0^+}\int_0^{\infty}e^{-xt}\varphi^{(r)}(t)\, dt=\int_0^{\infty}\varphi^{(r)}(t)\, dt=\infty,
$$
which contradicts \eqref{eq:ar}. Thus, $\varphi^{(r)}$ is identically zero and hence $\varphi$ must be a polynomial of degree $r-1$ with non-negative coefficients.

It follows from the discussion above that if $f$ is the Laplace transform of an absolutely monotonic function $\varphi(t)=\sum_{n=0}^{\infty}a_nt^n$  in which infinitely many of the numbers $a_n$ are strictly positive, then for all positive integers $k$ we have   
$$
\lim_{x\to 0^+}x^{k}f(x)=\infty,
$$
and since consequently all derivatives of $x^kf(x)$ are unbounded near zero we also have 
$$
\lim_{x\to 0^+}(-1)^{n+k}\left(x^{k}f(x)\right)^{n+k}=\infty
$$
for all $n,k\geq 0$.
\end{rem}
\begin{rem}
The generalized hypergeometric series
$$
\varphi(t)=\, _{1}F_{2}\big(a;b,c;t)=\sum_{k=0}^{\infty}\frac{(a)_k}{(b)_k(c)_kk!}t^k,\:\:a>0,b>0,c>0
$$
defines an absolutely monotonic function on $[0,\,\infty)$.
According to \cite[p. 115]{aar} (see also \cite{gasper}) its Laplace transform exists for all $x>0$ and it is given by the formula
\begin{equation}\label{f1}
f(x)=\int_{0}^{\infty}e^{-xt}\,\varphi(t)\,dt=\frac{1}{x}\;_{2}F_{2}\big(a,1;b,c;\frac{1}{x}\big)=\sum_{n=0}^{\infty}\frac{(a)_{n}}{(b)_{n}\,(c)_{n}}\,\frac{1}{x^{n+1}}.
\end{equation}
Therefore
$f$ has all the properties  of Theorem \ref{thm:main1}.
Notice that the case $a=b$, $c=1$ gives the example $H(x)=x^{-1}e^{1/x}$, mentioned in the first section. The case $a=b$, $c=\alpha+1$ gives
an example related to the modified Bessel function of the first kind. 
Formula \eqref{f1} has the following generalization
\begin{equation*}
\int_{0}^{\infty}e^{-xt}\,t^{\lambda-1}\,_{1}F_{2}\big(a;b,c;t)\,dt=\frac{\Gamma(\lambda)}{x^{\lambda}}\;_{2}F_{2}\big(a,\lambda;b,c;\frac{1}{x}\big)\,,
\end{equation*}
for any $\lambda>0$. Therefore the function $\tfrac{\Gamma(\lambda)}{x^{\lambda}}\;_{2}F_{2}\big(a,\lambda;b,c;\tfrac{1}{x}\big)$
has all the properties of Theorem \ref{thm:main4}.
\end{rem}

In order to give an example of a function satisfying the properties of Theorem 1.3 we describe a more general situation in the proposition below.

\begin{prop}
\label{prop:1}
Suppose that $f$  satisfies  the conditions in Theorem \ref{thm:main1}, let $\lambda>0$ and define $F_{\lambda}(x)=x^{-\lambda}f(x)$. Then: 
\begin{enumerate}
\item[(1)] If $\lambda$ is a positive integer then $F_{\lambda}$ satisfies again all the conditions in Theorem \ref{thm:main1}.
\item[(2)] If $\lambda$ is not an integer, let $k=[\lambda]$. Then $(-1)^{j}(x^{j}F_{\lambda}(x))^{(j)}$ is completely monotonic for $j\leq k+1$, so
$$
F_{\lambda}(x)=\int_0^{\infty}e^{-xt}p_{\lambda}(t)\, dt,
$$
with $t^kp_{\lambda}(t)\in \mathcal A_{k+1}(\sigma)$ for some $\sigma\in \mathcal A_0$.

Moreover, if $\lim_{x\to \infty}xf(x)>0$ then $(-1)^{k+2}(x^{k+2}F_{\lambda}(x))^{(k+2)}$ is {\em not} completely monotonic.
\end{enumerate}
\end{prop}
{\it Proof.} We have $f(x)=\sum_{n=0}^{\infty}a_nn!x^{-n-1}$ and hence (1) is obvious. To prove (2) notice that 
\begin{align*}
(-1)^{j}(x^{j}F_{\lambda}(x))^{(j)}&=(-1)^j\sum_{n=0}^{\infty}a_nn!(x^{j-n-1-\lambda})^{(j)}\\
&=\sum_{n=0}^{\infty}a_nn!(n+1+\lambda-j)\cdots(n+\lambda)x^{-n-1-\lambda}.
\end{align*}
For $j\leq k+1$ all coefficients in this series are non-negative, and this gives the complete monotonicity.
The representation then follows from Theorem \ref{thm:main2}. To obtain the last assertion we get from the formula above
$$
(-1)^{k+2}(x^{k+2}F_{\lambda}(x))^{(k+2)}=a_0(\lambda-k-1)(\lambda-k)\cdots\lambda x^{-\lambda-1}+o(x^{-\lambda-1})
$$
for $x\to \infty$. Here, $c_{\lambda}=(\lambda-k-1)(\lambda-k)\cdots\lambda$ is negative so that
$$
\lim_{x\to \infty}x^{\lambda+1}(-1)^{k+2}(x^{k+2}F_{\lambda}(x))^{(k+2)}=a_0c_{\lambda}<0,
$$
contradicting complete monotonicity.\hfill $\square$

\noindent
Stamatis Koumandos\\
Department of Mathematics and Statistics\\
The University of Cyprus \\ 
P.\ O.\ Box 20537\\
1678 Nicosia, Cyprus\\
{\em email}:\hspace{2mm}{\tt skoumand@ucy.ac.cy}

\vspace{0.5cm}

\noindent
Henrik Laurberg Pedersen\\
Department of Mathematical Sciences\\
University of Copenhagen \\
Universitetsparken 5\\
DK-2100, Denmark\\
{\em email}:\hspace{2mm}{\tt henrikp@math.ku.dk}


\begin{thebibliography}{99}
\bibitem{aar} G.E. Andrews, R. Askey, R. Roy, Special Functions, Cambridge Univ. Press, Cambridge, 1999.
\bibitem{boas} R.P.~Boas, Entire Functions. Academic Press, New York (1954).

\bibitem{gasper} G. Gasper, Positive integrals of Bessel functions, \emph{SIAM. J. Math. Anal.} 6(5) (1975) 868--881.

\bibitem{kp1} S.~Koumandos and H.L.~Pedersen, Completely
    monotonic functions of positive order and asymptotic expansions of the
    logarithm of Barnes double gamma
    function and Euler's gamma function, {\em J.\ Math.\ Anal.\
      Appl.}, {\bf 355} (2009), 33--40.

\bibitem{kp2}S.~Koumandos and H.L.~Pedersen, On Asymptotic Expansions of Generalized Stieltjes Functions, Comput.\ Methods Funct.\ Theory 15 (2015), 93-115 (DOI: 10.1007/s40315-014-0094-7).
 
\bibitem{SSV} R. L. Schilling, R.  Song and Z.  Vondra\v{c}ek, Bernstein functions. De Gruyter Studies in Mathematics \textbf{37}, De Gruyter, Berlin (2010).

\bibitem{Sok}  A. D. Sokal, Real-variables characterization of generalized Stieltjes functions,
\emph{Expo. Math.} \textbf{28} (2010), 179--185.

\bibitem{Wid1} D.V. Widder, The Stieltjes transform. \emph{Trans. Amer. Math. Soc.} \textbf{43}, (1938), 7--60.

\bibitem{Wid2} D.V. Widder,  The Laplace transform. Princeton University Press, Princeton (1946).
\end{thebibliography}
\end{document}